\documentclass[a4paper, 11pt,twosided]{article}
\usepackage[english]{babel}
\usepackage{amsmath, amssymb,theorem}
\usepackage[matrix,arrow]{xy}
\usepackage{epsfig}
\title{Some Hopf Algebras of Trees}
\author{Pepijn van der Laan \\ \small{vdlaan@math.uu.nl} \\ \small{Mathematisch Instituut,
Universiteit Utrecht}}
\date{Decemeber 9, 2001}
\newcommand{\note}{\textbf}

\newcommand{\NN}{\mathbb{N}}
\newcommand{\ZZ}{\mathbb{Z}}

\renewcommand{\phi}{\varphi}
\newcommand{\Del}{\Delta}

\newcommand{\eps}{\varepsilon}

\newcommand{\rarr}{\rightarrow}

\newcommand{\id}{\text{id}}

\newcommand{\Hom}{\text{Hom}}
\newcommand{\clcdot}{\cdot\ldots\cdot}

\newcommand{\be}{\begin{equation}}
\newcommand{\ee}{\end{equation}}
\newcommand{\uea}{universal enveloping algebra}
\newcommand{\hoa}{Hopf algebra}

\theoremheaderfont{\scshape}
\theorembodyfont{\itshape}
\theoremstyle{change}
\newtheorem{Tm}{Theorem}[section]
\newtheorem{Pp}[Tm]{Proposition}
\newtheorem{Lm}[Tm]{Lemma}
\newtheorem{Cr}[Tm]{Corollary}

\theorembodyfont{\upshape\rmfamily}
\newtheorem{Cv}[Tm]{Convention}

\newtheorem{Ex}[Tm]{Example}
\newtheorem{Rm}[Tm]{Remark}

\newenvironment{Pf}{{\scshape Proof}}{\hspace*{\fill}{\scshape QED}\vspace*{.4cm}}
\begin{document}

\maketitle

\section{Introduction}
\label{Sec:Intro}
In the literature several Hopf algebras that can be described in terms
of trees have been studied. This paper tries to answer the question
whether one can understand some of these Hopf algebras in terms of a
single mathematical construction. 

We recall the Hopf algebra of rooted 
trees as defined by Connes and Kreimer in \cite{ConKr:Hoa}. 
Apart from its physical relevance, it has a universal property in Hochschild 
cohomology. We generalize the operadic construction by Moerdijk
\cite{Moer:hoa} of this Hopf algebra to more general trees 
(with colored edges), and prove a universal property in coalgebra Hochschild 
cohomology. For a Hopf operad
$\textsf{P}$, the construction is based on the operad
$\textsf{P}[\lambda_n]$ obtained from \textsf{P} by adjoining a free
$n$-ary operation. 

For a specific multi-parameter family of Hopf algebras obtained in
this way, we give explicit formulas for the (a priori inductively defined) 
comultiplication and for the Lie bracket underlying the dual . 
Some special cases yield natural examples of pre-Lie and dendriform 
(and thus associative) algebras. We apply these results to construct
some known Hopf algebras of trees. Notably, the
Loday-Ronco Hopf algebra of planar binary trees \cite{LodRon:Trees}, and the
Brouder-Frabetti pruning Hopf algebra \cite{BrouFrab:Trees}.

The author plans to discuss the simplicial algebra
structure on the set of initial algebras with free $n$-ary operation
$P_n$ for an operad \textsf{P} with multiplication in future work, 
together with some relations of the $\textsf{P}[\lambda_n]$ to the
algebra Hochschild complex. 

The author is grateful to Ieke Moerdijk for suggesting the Hopf
\textsf{P}-algebras $P_n$  and their 
simplicial structure and for motivating as well as illuminating discussions,
and to Lo\"ic Foissy and Maria Ronco for pointing out errors in
previous versions of this manuscript.

%
%

\section{Preliminaries}
\label{Sec:Preliminaries}
\label{trees}
This section fixes notation on trees and describes some results of 
Connes and Kreimer\cite{ConKr:Hoa}, and Chapoton and Livernet
\cite{ChapLiv:Prelie} that motivated this paper.

\note{Rooted trees} $t$ are
isomorphism classes of finite partially ordered sets which
\begin{enumerate}
\item have a minimal element $r$ ($\forall x\neq r: r<x$); we call
$r$ the \note{root}, and 
\item satisfy the tree condition that 
$(y\neq z) \land (y<x) \land (z<x)$ implies $(y<z)\lor (z<y)$.
\end{enumerate}
The elements of a tree are called \note{vertices}. A pair of vertices
$v<w$ is called an \note{edge} if there is no vertex $x$ such that $v<x<w$.
The number of vertices of a tree $t$ is denoted by $|t|$.
A \note{path} from $x$ to $y$ in a tree is a sequence $(x_i)_i$ of
elements $x=x_n>x_{n-1}>\ldots>x_1>x_0=y$ of maximal length. We will
say that $x$ is above $y$ in a tree if there is a path from $x$ to $y$.
Thus we may depict a rooted tree as a finite directed graph, with one terminal
vertex, the root. A vertex is called a \note{leaf} if there is no
other vertex above it. A vertex is an \note{internal vertex} if it is
not a leaf. We draw trees the natural way, with the root downwards.
An automorphism of a tree $t$ is an automorphism of the partially 
ordered set.

A \note{forest} is a finite, partially ordered set satisfying only
property \textit{(ii)} of the above. A forest always is a disjoint 
union of trees, its connected components. A rooted tree can be pictured 
as a tree with unlabeled vertices and uncolored edges. A forest is a 
set of these.

In the sequel we need trees with colored edges. That is, there
is a function from the set of edges to a fixed set of colors.
Vertices can be labeled as well (i.e. there is a function from
the set of vertices to a fixed set of labels). Also trees
with a linear ordering on the incoming edges at each vertex are used.
A colored forest is a linearly ordered set of colored trees.
An automorphism of a colored tree is an automorphism of the underlying 
tree compatible with the coloring of the edges.

We consider two operations on trees: cutting and grafting.
Let $s$ and $t$ be trees and let $v\in t$ be a vertex of $t$. 
Define $s\circ_v t$ to be the tree that as a set is the disjoint
union of $s$ and $t$ endowed with the ordering obtained by
adding the edge $v<r_s$. This operation is called \note{grafting}
$s$ onto $v$. When using trees with colored edges, we write
$s\circ_v^i t$ for grafting $s$ onto $v$ by an edge of color $i$.
For trees with labeled vertices, grafting preserves the labels. The
operation of grafting $s$ on $t$ is visualized below.

\begin{figure}[!ht]
\begin{center}
\epsfig{file=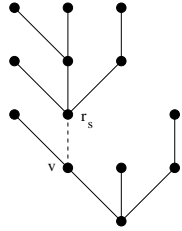,height=3cm,clip=}
\caption{grafting $s$ on $t$}
\end{center}
\end{figure}

A \note{cut} in a tree $t$ is a subset of the set of edges of $t$. 
A cut is \note{admissible} if for each leaf $m$ of $t$ the unique 
path $m>\ldots>x_2>x_1>r$ to the root contains at most one edge of the cut. 
Removing the edges in a cut $c$ yields a forest.
We denote by $R^c(t)$ the connected component of this forest 
containing the root and by $P^c(t)$ the complement of $R^c(t)$.
We denote the set of admissible cuts of $t$ by $C(t)$. For colored
and labeled trees, cutting preserves the colors of edges that
are not cut and labels of vertices. The left picture below shows an
admissible cut, whereas the cut on the right is not admissible.

\begin{figure}[!ht]
\begin{center}
\epsfig{file=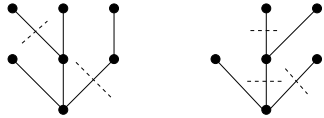,height=1.5cm,clip=}
\caption{an admissible cut (left) and a non-admissible cut (right)} 
\end{center}
\end{figure}

We describe the \note{\hoa\ $H_R$ of rooted trees} (defined in
Connes-Kreimer \cite{ConKr:Hoa}) as follows. Let $k$ be a field. 
The set of rooted trees and the empty set generate
a \hoa\ over $k$. As an algebra it is the polynomial algebra in formal
variables $t$, one for each rooted tree. Thus $H_R$ is spanned by forests.
Multiplication corresponds to taking the direct
union. The algebra $H_R$ is of course commutative. The unit is the
empty tree.

Comultiplication $\Delta$ is defined on rooted
trees $t$ and extended as a algebra homomorphism:
\[ 
\Delta(t)=\sum_{c\in C(t)}P^c(t)\otimes R^c(t), 
\]
where the sum is over admissible cuts.
The counit $\eps:H_R\rightarrow k$ takes value $1$ on the empty tree
and $0$ on other trees.

The \hoa\ $H_R$ is $\ZZ$-graded with respect to the number of vertices
in forests. The homogeneous elements of degree $m$ are products
$t_1\clcdot t_n$ of trees, such that $\sum |t_i| = m$. Both the product
and the coproduct preserve the grading.

Denote by $V^*$ the graded dual of a graded vector space $V$:
\[
V^* = \bigoplus_n (V_n)^*,
\]
the direct sum of the duals of the spaces of homogeneous
elements. A basis for $(H_R)_n^*$ is $\{D_{t_1\clcdot t_k}|\sum_i|t_i|=n\}$, 
the dual basis to the basis for $(H_R)_n$ given by products of trees. 
Of course, $H_R^*$ is a cocommutative Hopf algebra, with comultiplication
\[
\Delta(D_{t_1\clcdot t_k})
= \sum_{i=0}^k\sum_{\sigma\in S(i,k-i)}
D_{t_{\sigma^{-1}(1)},\ldots,t_{\sigma^{-1}(i)}}
\otimes D_{t_{\sigma^{-1}(i+1)},\ldots,t_{\sigma^{-1}(k)}},
\]
where $S(i,k-i)$ is the set of $(i,k-i)$-shuffles in $S_k$. We sum over
these to avoid repetition of terms in the right hand side. The primitive
elements are those dual to (single) rooted trees.

As proved by Connes and Kreimer \cite{ConKr:Hoa}, the Lie algebra which 
has $H_R^*$ as its \uea\ is the linear span of rooted trees with the 
Lie bracket
\[
[D_t,D_s] = \sum_{u} (n(t,s,u)-n(s,t,u)) D_{u},
\]
where $n(t,s,u)$ is the number of admissible cuts $c$ of $u$ such that
$P^c(u)=t$ and $R^c(u)=s$.
Chapoton and Livernet \cite{ChapLiv:Prelie} show that this is the Lie
algebra associated to the free pre-Lie algebra on one generator (cf. 
section \ref{Sec:PreLie} for details).

%
%

\section{Operads}
\label{Sec:Operads}
\begin{Cv}
For the rest of this paper, we restrict to the category of vector
spaces over a fixed field $k$, but the general theory carries over
to any symmetric monoidal category with countable coproducts and quotients 
of actions by finite groups on objects
(cf. Moerdijk \cite{Moer:hoa}).
\end{Cv}

Let $kS_n$ denote the group algebra of the permutation group $S_n$ on $n$ objects.
An \note{operad} \textsf{P} in this paper will mean an operad with unit. 
Thus an operad consists of a collection of right $kS_n$-modules
$\textsf{P}(n)$, together with an associative, $kS_n$-equivariant composition
\[
\gamma:\textsf{P}(n)\otimes\textsf{P}(m_1)\otimes\ldots\otimes\textsf{P}(m_n)
\longrightarrow \textsf{P}(m_1+\ldots + m_n),
\]
such that there exists an identity $\id:k\rarr\textsf{P}(1)$ with the obvious 
property and a unit $u:k\rarr\textsf{P}(0)$, which need not be an isomorphism.
The existence of the unit map is the only change with respect to
the usual definition of an operad (consult Kriz and May
\cite{KrizMay:Opd}, or Ginzburg and Kapranov \cite{GinKap:Koszul},
or Getzler and Jones \cite{GetzJon:Opd}). Following Getzler and Jones 
\cite{GetzJon:Opd} we define an \note{operad with multiplication} as an operad 
together with an associative element $\mu\in \textsf{P}(2)$.

The category of \note{pointed vector spaces} has as objects vector
spaces $V$ with a base point $u:k\rarr V$. Morphisms are base
point preserving $k$-linear maps. We use notation $1:=u(1)$. There
are some canonical functors: The forgetful functor to
vector spaces will be used implicitly. The free 
associative algebra functor $T$ on a pointed vector space $V$ is left
adjoint to the forgetful functor from unital associative algebras to
pointed vector spaces. 

For any pointed vector space $V$, there is an operad
$\textsf{End}_V$ such that $\textsf{End}_V(n)=\Hom_k(V^{\otimes
n},V)$. The pointed structure of $V$ is only used to
to define $u:k\rightarrow \textsf{End}(V)(0)$, the elements of 
$\Hom_k(V^{\otimes n},V)$ need not preserve the base point.
A \textsf{P}\note{-algebra} structure on $V$ is a map of operads
(i.e. preserving relevant structure) from \textsf{P} to
$\textsf{End}_V$. Equivalently, a \textsf{P}-algebra is a vector space $V$
together with linear maps
\[
\gamma_V:\textsf{P}(n)\otimes_{kS_n}( V^{\otimes m_1}\otimes\ldots\otimes V^{\otimes m_n})
\longrightarrow V^{\otimes m_1+\ldots m_n},
\]
compatible with composition in \textsf{P} and unit in the natural sense.
\begin{Ex}
Let \textsf{k} be the operad, with $\textsf{k}(0) = k = \textsf{k}(1)$
and $\textsf{k}(n)=0$, otherwise (with $u=\id$). Algebras for
\textsf{k} are pointed vector spaces. 
The operad \textsf{k} is the initial operad in our setting.

We denote by \textsf{Com} the operad with as algebras unital
commutative algebras. This operad satisfies $\textsf{Com}(n)=k$ for
all $n$ (where $u:k\rightarrow k$ is the identity map).

Likewise \textsf{Ass} is the operad satisfying $\textsf{Ass}(n) =
kS_n$, the group algebra of $S_n$, as a right $S_n$-module. The
\textsf{Ass}-algebras are unital associative algebras.
An operad with multiplication is an operad \textsf{P} together with a
map of operads from \textsf{Ass} to \textsf{P}.
\end{Ex}
A map of operads
$\phi:\textsf{P}\rightarrow\textsf{Q}$ induces an obvious functor
$\phi^*:\textsf{Q}\text{-Alg}\rightarrow\textsf{P}\text{-Alg}$.
This functor has a left adjoint
$\phi_!:\textsf{P}\text{-Alg}\rightarrow\textsf{Q}\text{-Alg}$.
Let \textsf{P} be any operad and let $i:\textsf{k}\rightarrow \textsf{P}$ 
be the unique inclusion of operads. Then $i_!$ is the unitary free 
\textsf{P}-algebra functor.

Let \textsf{P} be an operad. We can form the operad
$\textsf{P}[\lambda_n]$ by adjoining a free $n$-ary operation to
$\textsf{P}(n)$. Algebras for
$\textsf{P}[\lambda_n]$ are just \textsf{P}-algebras $A$ endowed
with a linear map $\alpha: A^{\otimes n}\rarr A$.
The reader familiar with the collections (cf. Getzler and 
Jones \cite{GetzJon:Opd}, or Ginzburg and Kapranov
\cite{GinKap:Koszul}), will recognize $\textsf{P}[\lambda_n]$ as the
coproduct of operads
\[
\textsf{P}[\lambda_n] = \textsf{P}\oplus_{\textsf{k}}FE_n
\]
of $\textsf{P}$ and the free operad $FE_n$ on the collection $E_n$ 
defined by $E_n(n) = kS_n$ and $E_n(m)=0$ for $m\neq n$.

Denote the initial \textsf{P}$[\lambda_n]$-algebra
\textsf{P}$[\lambda_n](0)$ by $P_n$.
For any \textsf{P}$[\lambda_n]$-algebra $(A,\alpha)$, there is by
definition an unique \textsf{P}-algebra morphism $\gamma$ such that
the following diagram commutes.
\[
\xymatrix{ P_n^{\otimes n} \ar[r]^{\lambda_n}\ar[d]_{\gamma^{\otimes n}} & P_n \ar[d]^{\gamma}\\
A^{\otimes n} \ar[r]^{\alpha} & A}
\]

The operad $\textsf{P}[\lambda_n]$ can be described in terms of planar
trees.
A planar tree is understood to have a linear ordering on the incoming 
edges at each vertex and might have external edges. There is a linear
ordering on all external edges. An element of
$\textsf{P}[\lambda_n](m)$ can be represented by a tree with $m$
external edges
and each vertex $v$ with $n$ incoming edges labeled either by
an element of $\textsf{P}(n)$ or by $\lambda_n$.
We identify two such representing trees if they can be reduced to the same 
tree using edge contractions implied by composition in \textsf{P}, 
equivariance with respect to the $kS_{n}$-action at each internal
vertex $v$ and $kS_m$-equivariance at the external edges. The
initial $\textsf{P}[\lambda_n]$-algebra $\textsf{P}[\lambda_n](0)$ 
is described by such trees without external edges.

\begin{Ex}
Let \textsf{P = Com}. The initial algebra $C_n = \textsf{Com}[\lambda_n](0)$
can be described as the free commutative algebra on trees with edges colored by
$\{1,\ldots,n\}$ 
(and no ordering on the edges). A bijection $T$ is
given by induction on the
number of applications of $\lambda$ as follows.
$T(1) = 1=\emptyset$, and $T(\lambda(1,\ldots,1)) = r$ (the one vertex
tree), and $T(\lambda(x_1,\ldots,x_n)$ is the 
tree obtained from the forest $x_1\clcdot x_n$ by adjoining a new
root and connecting the roots of each tree in $x_i$ to the new root by an edge 
of color $i$.
Figure \ref{Fig:ComAss} shows (twice) the tree
$T(\lambda(x,\lambda(xx,1)\lambda(x,1)))$
in $C_2$, where $x=\lambda(1,1)$.

\begin{figure}[!ht]
\begin{center}
\epsfig{file=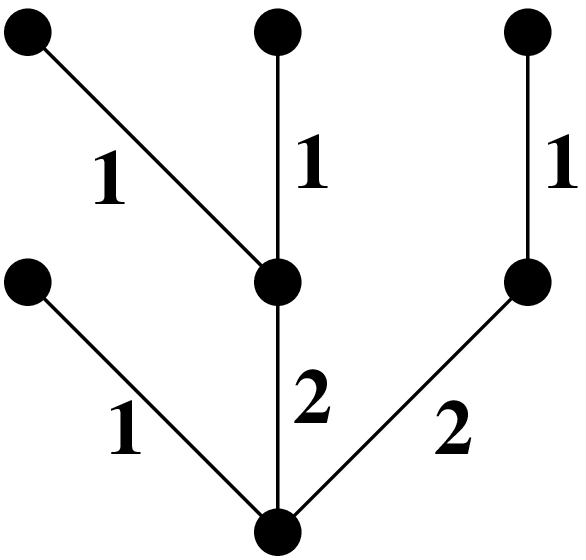,height=1.5cm,clip=}
\hspace*{1.5cm}
\epsfig{file=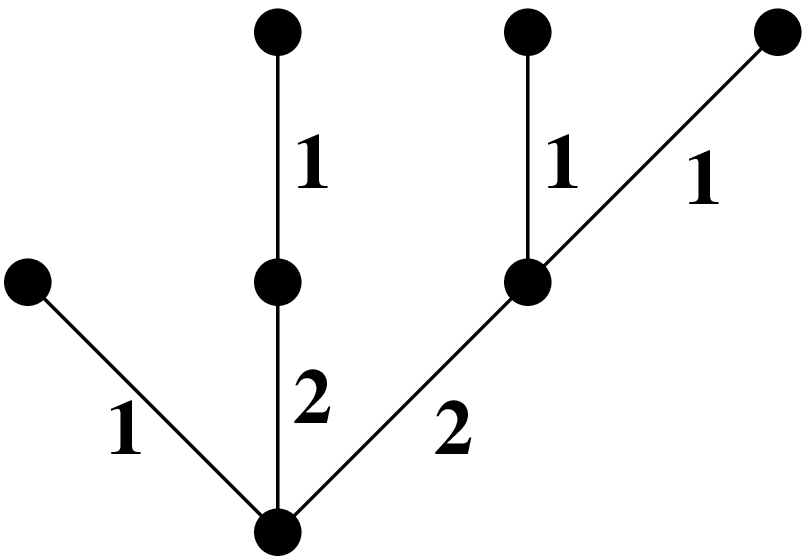,height=1.5cm,clip=}
\end{center}
\caption{trees not equal in $A_2$, but equal in
$C_2$\label{Fig:ComAss}}
\end{figure}

\end{Ex}
\begin{Ex}
Let $\textsf{P = Ass}$. The initial algebra $A_n = \textsf{Ass}[\lambda_n](0)$
can be described as the free associative algebra on trees with edges colored 
by $\{1,\ldots,n\}$ and at each vertex the incoming edges of each color
endowed with a linear ordering.
Equivalently, we could say, a linear ordering on all incoming edges extending 
the ordering on the colors. A bijection is given by the map $T$
with the same inductive definition as in the previous example, 
but using associative multiplication. In the associative case, 
the two trees in figure \ref{Fig:ComAss} are not identified.
The tree to the right in figure \ref{Fig:ComAss} is 
$T(\lambda(x,\lambda(x,1)\lambda(xx,1)))$ in $A_2$, while the tree to
the left is $T(\lambda(x,\lambda(xx,1)\lambda(x,1)))$ (where $x=\lambda(1,1)$).
\end{Ex}

\begin{Ex}
Let $\textsf{P = k}$. Then $k_n = \textsf{k}[\lambda_n](0)$ can be
identified with the free vector space on trees (not forests, but
possibly empty) with edges labeled by $\{1,\ldots,n\}$ and at most one
incoming edge of each color at each vertex. For the relation to
$n$-ary trees, see section \ref{Sec:Pruning}.
\end{Ex}

%
%

\section{Hopf Operads}\label{Sec:HopfOperads}

A \note{Hopf operad} (cf. Getzler and Jones \cite{GetzJon:Opd}, and
also Moerdijk \cite{Moer:hoa}) is an operad
in the category of coalgebras over $k$. Thus, a \note{Hopf operad}
\textsf{P} is a ($k$ linear) operad together with a morphisms
\[
k\overset{\eps}{\longleftarrow} \textsf{P}(n)
\overset{\Delta}{\longrightarrow}\textsf{P}(n)\otimes\textsf{P}(n),
\]
such that $\Delta$ and $\eps$ satisfy the usual axioms of a
coalgebra: $\Delta$ is coassociative and $\eps$ is a counit for
$\Delta$ (cf. Sweedler \cite{Sweedler:Hoa}), and such that the
composition $\gamma$ of the operad is a coalgebra morphism. Moreover, 
$\Delta$ should be compatible with the $kS_n$-action, where 
$\textsf{P}(n)\otimes\textsf{P}(n)$ is a $kS_n$- module via the 
diagonal coproduct. For any Hopf operad \textsf{P}, the tensor product of two
\textsf{P}-algebras is a \textsf{P}-algebra again. 

A Hopf
\text{P}-algebra $A$ is a \textsf{P}-algebra $A$ in the category
of (counital!) coalgebras. Note that we can not use the description of
algebras in terms of the endomorphism operad in the category of
coalgebras, since coalgebra homomorphisms do not form a linear
space. In this generality it is not natural to consider antipodes. A
Hopf \textsf{Ass}-algebra is just a bialgebra.

Let $\phi:\textsf{Q}\rarr \textsf{P}$ be a map of
Hopf operads. The map $\phi$ induces functors
\[
\phi^*:\textsf{P}\text{-Alg}\rightarrow \textsf{Q}\text{-Alg}
\qquad \text{and} \qquad 
\bar{\phi}^*:\textsf{P}\text{-HopfAlg}\rightarrow
\textsf{Q}\text{-HopfAlg}.
\]
The map $\phi^*$ has a left adjoint $\phi_!$ (we work 
with $k$-vector spaces). Note that $\phi^*(k)=k$ and for 
\textsf{P}-algebras $A$ and $B$ we have $\phi^*(A\otimes B) = 
\phi^*(A)\otimes \phi^*(B)$ (the map $\phi$ is compatible with $\Delta$).
Using this observation, conclude that the adjunction induces algebra maps 
$\phi_!(k)=\phi_!\phi_*(k)\rarr k$ and
$\phi_!(A\otimes B)\rarr \phi_!(\phi^*\phi_!A\otimes\phi^*\phi_!B)
\rightarrow\phi_!(A)\otimes\phi_!(B)$. These maps
serve to show that $\phi_!$ lifts to a left adjoint
\[
\bar{\phi}_!:\textsf{Q}\text{-HopfAlg}\longrightarrow
\textsf{P}\text{-HopfAlg}
\]
of $\bar{\phi}^*$.

Let $\textsf{P}$ be a Hopf operad. The aim of this section is to
state a general result on Hopf operad structures on $\textsf{P}[\lambda_n]$.
We use the notations $\lambda=\lambda_n$ and
$P_n=\textsf{P}[\lambda_n](0)$ when no confusion can arise.

Let $(A,\alpha)$ be a $\textsf{P}[\lambda_n]$-algebra and 
let $\sigma_1,\sigma_2: A^{\otimes n}\rarr A$ be a pair of
linear maps and define for each such pair a map
$(\sigma_1,\sigma_2):(A\otimes A)^{\otimes n}\rarr A\otimes
A$ by
\[
(\sigma_1,\sigma_2) = (\sigma_1 \otimes \alpha_n
+ \alpha_n \otimes \sigma_2)\circ\tau,
\]
where $\tau$ is the $n$-fold
twist map which identifies $(A\otimes A)^{\otimes n}$ with
$A^{\otimes n}\otimes A^{\otimes n}$ by the symmetry of 
the tensor product. For vector spaces $V_1,\ldots,V_{2n}$,
the map $\tau$ is the natural map 
\[
\tau:V_1\otimes\ldots\otimes V_{2n}\mapsto V_1\otimes\ldots
V_{2n-1}\otimes V_2\otimes V_4\ldots\otimes V_{2n}
\]
induced by the symmetry (in the dg case this involves some signs). 

We now list some results for linear maps
$\sigma_i$ that define a Hopf algebra structure on $P_n$ but do not
extend to $\textsf{P}[\lambda_n]$. The results in this section are
obtained along the same lines as the results in Moerdijk \cite{Moer:hoa}.
Our approach is more general in the sense that Moerdijk considered
only the case $n=1$.

Let $(P_n,\lambda_n)$ be the initial
$\textsf{P}[\lambda_n]$-algebra. For any pair of $n$-ary linear maps 
$\sigma_1,\sigma_2$ there is a unique $\textsf{P}[\lambda_n]$-algebra morphism
$(P_n,\lambda_n)\rightarrow (P_n\otimes P_n,(\sigma_1,\sigma_2))$. That
is, a unique $P$-algebra morphism such that
\[
\xymatrix{ P_n^{\otimes n} \ar[r]^{\lambda}\ar[d]^{\Delta^{\otimes n}} 
& P_n\ar[d]^{\Delta} \\
(P_n\otimes P_n)^{\otimes n}\ar[r]^{\ (\sigma_1,\sigma_2)}
& P_n\otimes P_n}
\]
commutes. Define $\eps:P_n\rightarrow k$ as the unique
$\textsf{P}[\lambda]$-algebra morphism $(P_n,\lambda_n)\rightarrow (k,0)$,
which extends $\eps:P(0)\rightarrow k$.
The following proposition is a generalization of Moerdijks construction 
\cite{Moer:hoa}. The proof is completely analogous. 

\begin{Pp} \label{Tm:genmoerdijk}
Let $n \in \NN$, $\lambda_n$ and $P_n$ be defined as above. Let
$\sigma_i:P_n^{\otimes n}\rarr P_n$ for $i=1,2$ be linear maps. 
If $\sigma_i$ satisfies
\[
\begin{split}
\epsilon \circ \sigma_i &= \epsilon^{\otimes n},\quad\text{and}\\
(\sigma_i\otimes\sigma_i)\circ\tau\circ\Delta^{\otimes n} &=
\Delta\circ\sigma_i;
\end{split}
\]
then there exists a unique Hopf-$\textsf{P}$ algebra structure on $P_n$
extending the Hopf $P$-algebra structure on $P(0)$ such that 
$\Delta\circ\lambda = (\sigma_1,\sigma_2) \circ \tau \circ
\Delta^{\otimes n}$ and $\eps\circ\lambda = 0$.
\end{Pp}

In the following examples $|t|$ denotes the number of applications
of $\lambda_n$ in an element $t\in\textsf{P}[\lambda_n](m)$, for any
$m$. In the case $\textsf{P} = \textsf{Ass}$, or $\textsf{P} =
\textsf{Com}$, or $\textsf{P} = \textsf{k}$ this corresponds to
counting the number of vertices in a tree. 

\begin{Ex}
Let \textsf{P} be a Hopf operad. We consider the case $n=1$.

i) For any $q_1,q_2\in k$, the endomorphisms $\sigma_i(t)=q_i^{|t|}t$
define a coassociative Hopf $\textsf{P}$-algebra structure on 
$P_1$ (cf. 
Moerdijk \cite{Moer:hoa} for the initial algebra).

ii) Suppose that \textsf{P} is a Hopf operad with multiplication. For any
$q_1,q_2\in k$, the endomorphisms
$\sigma_i(t)=q_i^{|t|}r^{|t|}$ define a Hopf $\textsf{P}$-algebra 
structure on $P_1$. (Here $r$ is $\lambda(1)$,
and powers are with respect to the multiplication in
\textsf{P}.)
\end{Ex}

\begin{Ex}\label{Ex:deform}
Let $\textsf{P}$ be a Hopf operad with multiplication. 
For any choice of $q_{i\,j}\in k$ for all $j\leq n$, the
maps
\[
\sigma_i(t_1,\ldots,t_n)=q_{i\,1}^{|t_1|}\clcdot q_{i\,n}^{|t_n|}
t_1\clcdot t_n
\]
(where $i=1,2$) give a Hopf $\textsf{P}$-algebra structure on $P_n$.
\end{Ex}

\section{Cohomology}
\label{Sec:Cohomology}
In the previous section we constructed Hopf operad structures on 
$\textsf{P}[\lambda_n]$. 
This section aims to describe the relation of these Hopf operads
and Hochschid cohomology for coalgebras.

Let $A$ be a Hopf \textsf{Ass}-algebra, 
and let $C^*_{(p)}$ be the graded vector space, which in degree $q$ is 
$\Hom_k(A^{\otimes p},A^{\otimes q})$.
For $\phi\in\Hom_k(A^{\otimes p},A^{\otimes q})$ 
we define a differential by the formula
\[
d\phi = (\mu^{(p)}\otimes\phi)\circ\tau\circ\Del^{\otimes p} +
\sum_{i=1}^q(-1)^i\Del_{(i)}\circ\phi
+(-1)^{q+1}(\phi\otimes\mu^{(p)})\circ\tau\circ\Del^{\otimes p},
\]
where $\Del_{(i)}(x_1\otimes\ldots\otimes x_q) = (x_1\otimes
\ldots\otimes\Del(x_i)\otimes\ldots\otimes x_q)$. This is the 
coalgebra-Hochschild complex with respect to
the left and right coaction of $A$ on $A^{\otimes p}$ given by
\[
(\mu^{(p)}\otimes\id)\circ\tau\circ\Delta: A^{\otimes
p}\longrightarrow A\otimes A^{\otimes p}
\]
(or with $(\id\otimes\mu^{(p)}$), where $\mu^{(p)}$ is the unique
$p$-fold application of $\mu$. This complex (for $p,q>0$) is the $p$-th 
column in the bicomplex used by Lazarev and Movshev
\cite{LazMov:hoa} to compute what they call the cohomology of a Hop
algebra. 

Let \textsf{P} be a Hopf operad, and let $A$ be a Hopf \textsf{P}-algebra, 
let $p\geq 1$ and $\sigma_1,\sigma_2:A^{\otimes p}\rarr A$ 
be coalgebra morphisms. Then
\[
\begin{split}
(\sigma_1\otimes\id)\circ \tau\circ\Delta^{\otimes p}:& A^p\rarr
A\otimes A^{\otimes p} \qquad\text{ and }\\
(\id\otimes\sigma_2)\circ \tau\circ\Delta^{\otimes p}:& A^p\rarr
A^{\otimes p}\otimes A
\end{split}
\]
define a left and a right coaction of $A$ on $A^p$ and we can
define the complex $C^q_{\sigma_1\sigma_2}$ with the boundary
$d_{\sigma_1\sigma_2}$. We can write the
differential explicitly as
\[
d\phi = (\sigma_1\otimes\phi)\circ\tau\circ\Del^p +
\sum_{i=1}^q(-1)^i\Del_{(i)}\circ\phi
+(-1)^{q+1}(\phi\otimes\sigma_2)\circ\tau\circ\Del^p.
\]
This is the Hochschild boundary with respect to these coactions. The 
cohomology of this complex will be denoted $H^{*}_{\sigma_1\sigma_2}(A)$.

Let \textsf{P} be a Hopf operad. A \note{natural $n$-twisting function} 
$\phi$ is a map from Hopf \textsf{P}-algebras $B$ to 
coalgebra maps $\phi^{(B)}:B^{\otimes n}\rarr B$, such that the map 
$\phi$ commutes with augmented \textsf{P}-algebra morphisms $f:A\rarr B$ 
(i.e. $f \circ\sigma_i^{(A)} = \sigma_i^{(B)} \circ f^{\otimes n}$ for
$i=1,2$).

If $\sigma_1$ and $\sigma_2$ are natural $p$-twisting functions
such that $\sigma^{(B)}_i$ satisfies the conditions from
\ref{Tm:genmoerdijk} for any Hopf \textsf{P} algebra $B$; then
$H^{*}_{\sigma_1\,\sigma_2}(B)$ is defined for any Hopf-\textsf{P}
algebra $B$, and natural in $B$.

Let $(B,\beta)$ be a Hopf \textsf{P}-algebra. Then $\beta$ is a 1-cocycle 
in the $(\sigma_1,\sigma_2)$-complex iff
\[
(\sigma_1\otimes\beta)\Delta - \Delta\beta +(\beta\otimes \sigma_1)\Delta =0.
\]
\begin{Ex}
Let $\textsf{P}$ be a Hopf operad. Then $\sigma=p\in\textsf{P}(n)$ 
defines a natural $n$-twisting function. More generally, if $\sigma=
(\eps^{\otimes k} \otimes p) \tau$ for some $p\in \textsf{P}(n-k)$ and
some $\tau\in S_n$, then $\sigma$ defines a natural $n$-twisting function.
\end{Ex}

We now characterize the  universal property of the Hopf
$P$-algebra $P_n$. Again this is a direct generalization of Moerdijk
\cite{Moer:hoa}. 

\begin{Tm}
Let $(\sigma_1,\sigma_2)$ be a pair of natural twisting functions of
Hopf $\textsf{P}$-algebras, satisfying the conditions of Theorem
\ref{Tm:genmoerdijk} and endow $P_n$ with the induced Hopf $\textsf{P}$-algebra
structure.
For every Hopf $\textsf{P}$-algebra $B$, there is a one-one correspondence
between 
\begin{enumerate}
\item Hopf $\textsf{P}$-algebra maps $c_\beta:P_n\rightarrow B$ that lift to a
map $(P_n,\lambda)\rightarrow (B,\beta)$ of $\textsf{P}[\lambda_n]$-algebras and 
\item cocycles $\beta$ in the $(\sigma_1,\sigma_2)$-twisted complex of
$B$.
\end{enumerate}
\end{Tm}

\begin{Pf}
The previous paragraph shows that $\lambda$ is a cocycle in the
$(\sigma_1,\sigma_2)$-complex of $P_n$ with respect to the
$(\sigma_1,\sigma_2)$-coproduct on $P_n$.

First suppose that $c_\beta$ has such a lift
$c_\beta:(P_n,\lambda)\rightarrow (B,\beta)$ (which is obviously
unique). Since $\sigma_1$ and
$\sigma_2$ are natural $n$-twisting functions and $c_\beta$ is a map
of augmented $\textsf{P}$-algebras, we have the following
commutative diagram of $\textsf{P}[\lambda_n]$-algebras:
\[
\xymatrix {(P_n,\lambda)\ar[rr]^{\Delta}\ar[d]_{c}& &
(P_n\otimes P_n,(\sigma_1,\sigma_2))\ar[d]^{c\otimes c} \\
(B,\beta)\ar[rr]_{\Delta}& & (B\otimes
B,(\sigma_1,\sigma_2)).}
\]
This shows that $\Delta\circ \beta = (\sigma_1,\sigma_2)\circ \Delta$,
which is the cocycle relation for $\beta$ in the
$(\sigma_1,\sigma_2)$-twisted complex.

Let $\beta$ be a cocycle in the $(\sigma_1,\sigma_2)$-twisted
complex of $B$. Then (by the universal property of $P_n$) there is a
unique $\textsf{P}[\lambda_n]$-algebra map $c:(P_n,\lambda)\rightarrow (B,\beta)$. 
We are done as soon as we have verified that the diagram above is
a commutative diagram of $\textsf{P}[\lambda_n]$-algebras. 
To show commutativity of the diagram, note that the universal
property of $P_n$ shows that $\eps\circ c = \eps$ is the unique
morphism of $(P_n,\lambda)$ to $(k,0)$. The intertwining relation
for natural $n$-twisting functions with respect to augmented
P algebra morphisms now implies that $c\otimes c$ is a
morphism of $\textsf{P}[\lambda_n]$-algebras. Since both $c$ and $\Delta$ are
maps of $\textsf{P}[\lambda_n]$-algebras, the universal
property of $P_n$ implies that both composites in the diagram are
the unique $\textsf{P}[\lambda_n]$-algebra morphism from
$(P_n,\lambda)$ to $(B\otimes B,(\sigma_1,\sigma_2))$.
\end{Pf}

The result above tells us that $\textsf{P}[\lambda_n]$ has the following
universal property in twisted coalgebra Hochschild cohomology, which
we state once more in $\textsf{P}[\lambda_n]$-free language.

\begin{Cr}\label{Cr:genmoerdijk}
Let $(B,\beta)$ be a pair of a Hopf $\textsf{P}$-algebra together with a
cocycle in the $(\sigma_1,\sigma_2)$-twisted complex.
Then there is a unique Hopf $\textsf{P}$-algebra morphism $c:P_n\rarr B$ 
such that the following diagram commutes:
\[
\xymatrix{ P_n^{\otimes n}\ar[d]^{c^{\otimes n}} \ar[r]^{\lambda} & P_n \ar[d]^{c} \\
B^{\otimes n} \ar[r]_{\beta} & B.}
\]
\end{Cr}

This section concludes with some functorial properties of the
construction. The proof of the following proposition is completely analogous 
the the proof for the case $n=1$ in Moerdijk \cite{Moer:hoa}.

\begin{Pp}\label{Pp:Functor}
Let $\phi: Q \rarr P$ be a map of Hopf operads, and 
let $\sigma_1,\sigma_2$ satisfy the conditions
of \ref{Tm:genmoerdijk} on $Q_n$ such that $\sigma_1\, \sigma_2$
extend to natural $n$-twisting functions of
$Q$-algebras; then
\begin{enumerate}
\item
The $\sigma_i$ induce natural twisting functions $\sigma_i$ of
$\textsf{P}$-algebras, satisfying the conditions of
\ref{Tm:genmoerdijk} on $H$.
\item
The unique map of $Q[\lambda_n]$-algebras $j_0:Q_n\rarr
\phi^*(P_n)$ is a map of Hopf $Q$-algebras.
\item
The unique map of $\textsf{P}[\lambda_n]$-algebras
$j:\phi_!(Q_n)\rightarrow P_n$ is a map of Hopf $\textsf{P}$-algebras.
\end{enumerate}
\end{Pp}

\begin{Ex}\label{Ex:Twistfunc}
Let $\sigma= (\eps^{\otimes m}\otimes p) \tau$ for some $m\leq n$, some
$p\in \textsf{P}(n-m)$ and some $\tau\in S_n$. If $p$ is group-like a
group-like element in the sense that $\Delta(p) = p\otimes p$, then
$\sigma$ satisfies the
conditions of theorems \ref{Tm:genmoerdijk} and \ref{Cr:genmoerdijk}.
In this case maps of operads $\psi:\textsf{Q}\rightarrow
\textsf{P}$ and $\phi:\textsf{P}\rightarrow
\textsf{Q}$ induce functors $\textsf{P}\text{-HopfAlg}\rightarrow
\textsf{Q}\text{-HopfAlg}$ that map the natural
twisting function $\sigma$ to a natural twisting functions
$\bar{\psi}^*\sigma$ and $\bar{\phi}_!\sigma$ of Hopf 
\textsf{Q}-algebras. 
\end{Ex}

%
%
%
%

\section{A Deformation}
\label{Sec:Deformation}
Moerdijk \cite{Moer:hoa} gives an inductive formula for a 2-parameter
deformation of $\Delta$ on $P_1$ for any Hopf operad \textsf{P}. 
Regarding $P_n$, 
we have seen in example \ref{Ex:deform} that there is a $2n$-parameter
deformation of the 
coproduct given by $\sigma_1=\sigma_2=u\circ\eps^{\otimes n}$ (provided 
that \textsf{P} is an Hopf operad with multiplication).
This section obtains an explicit formula for
the $2n$-parameter deformation of the coproduct in $C_n$ and the
Lie algebra structure underlying the dual. 

\begin{Cv}\label{Cv:Deform}
In this section we fix \textsf{P = Com}, and $n\in \NN$.
We use notation $C_n = \textsf{P = Com}[\lambda_n](0)$.
For $i=1,2$ and $1\leq j \leq n$, let $q_{ij}\in k$, and define
\[
\sigma_i(t_1,\ldots t_n) = t_1\clcdot t_n\cdot \prod_j
q_{ij}^{|t_j|}.
\]
According to example \ref{Ex:deform}, these $\sigma_i$ make $C_n$ a
commutative Hopf algebra.
\end{Cv}

A subforest $s$ of a rooted tree $t$ is a subset of the partially ordered set 
$t$ with the induced partial ordering. 
For colored trees, the color of the edge connecting $v>w$ in $s$ is the 
color of edge connecting $w$ to its direct predecessor in the unique 
path from $v$ to $w$ in $t$. For $v\in s$ we denote by $p_k(v,s,t)$ the number of 
edges of color $k$ in the path in $t$ from $v$ to the root of $t$ that have their 
lower vertex in $s^c$. For forests $t$ we define $p_k(v,s,t)$ as
$p_k(v,s\cap t',t')$, where $t'$ is the connected component of $t$
containing $v$.
There is an easy but useful lemma on the calculus of the $p_k$.
\begin{Lm}
Let $t$ and $s$ be subforests of a forest $u$. Let $v\in s$ and set 
$t' = t\cup v$, $s' = s\cap t'$, $t''= t^c\cup v$ and $s'' = s\cap t''$. Then 
\[
p_k(v,s,u) = p_k(v,s',t') + p_k(v,s'',t''),
\]
where $t',\ t'',\ s'$ and $s''$ are interpreted as subforests of $u$. 
\end{Lm}
\begin{Pf}
The lemma follows at once when we observe that a vertex in the path from 
$v$ to the root in $u$ that is not in $s$ is either in $t'$ or in $t''$.
\end{Pf}
Define
\[
q(s,t):=
\prod_j q_{1\,j}^{\sum_{v\in s}p_j(v,s,t)}\cdot
\prod_jq_{2\,j}^{\sum_{v\in s^c}p_j(v,s^c,t)}.
\]
Intuitively, $q(s,t)$ counts for $v\in s$ the number of edges of color $j$
are in the path from $v$ to the root that have their lower vertex in
$s^c$ and adds a factor $q_{1j}$ for each of these, and $q(s,t)$
counts for $v\in s^c$ the number of edges of color $j$
are in the path from $v$ to the root that have their lower vertex in
$s$ and adds a factor $q_{2,j}$ for each of these.

\begin{Tm}\label{Tm:Coprod}
Under the assumptions of \ref{Cv:Deform}, we have for a forest $t\in
C_n$ the formula
\[
\Delta(t) = \sum_{s\subset t} 
q(s,t) \, s\otimes s^c ,
\]
where the sum is over all subforests $s$ of $t$.
\end{Tm}
\begin{Pf}
We use induction with respect to the number of applications of $\lambda$. 
The formula is trivial for the empty tree.

Let $t=\lambda(x_1,\ldots,x_n)$ be a tree and suppose that the formula holds for all 
forests with less then $|t|$ vertices. Subforests of $t$ are either of the form
$ s = \cup_i s_i$, a (disjoint) union of subforests of the $x_i$, or of the form
$s = r \cup(\cup_i s_i)$, a (disjoint) union of subforests of the $x_i$ together with the root.
By definition,
\[
\begin{split}
\Delta(t) &= \sum_{s_i\subset x_i} s_1\clcdot s_n \otimes \lambda(s_1^c,\ldots,s_n^c)
\prod_i q_{1i}^{|s_i|} q(s_i,x_i)\\
&+ \sum_{s_i\subset x_i} \lambda(s_1,\ldots,s_n) \otimes s_1^c \clcdot s_n^c
\prod_i q_{2i}^{|s^c_i|} q(s_i,x_i).
\end{split}
\]
But by the lemma above,
\[
\prod q_{1i}^{|s_i|} q(s_i,x_i) = 
\prod_j q_{1\,j}^{\sum_{v\in s}p_j(v,s,t)}\cdot
\prod_jq_{2\,j}^{\sum_{v\in s^c}p_j(v,s^c,t)}
\]
for $s=\cup_i s_i = s_1\clcdot s_n$ and $s^c = r\cup(\cup_is_i^c)= \lambda(s_1^c,\ldots,s_n^c)$; and 
\[
\prod q_{2i}^{|s^c_i|} q(s_i,x_i) = 
\prod_j q_{1\,j}^{\sum_{v\in s}p_j(v,s,t)}\cdot
\prod_jq_{2\,j}^{\sum_{v\in s^c}p_j(v,s^c,t)}
\]
for $s= r\cup(\cup_i s_i) = \lambda(s_1,\ldots,s_n)$ and $s^c = \cup_is_i=s_1^c\clcdot s_n^c$. 
Putting these together yields the formula.
\end{Pf}

An \note{antipode} for a Hopf \textsf{Ass}-algebra $A$ is
an inverse for the identity map $id$ with respect to the convolution
product
\[
f*g = \mu\circ(f\otimes g)\circ \Delta
\]
on the $k$-linear maps $\Hom_k(A,A)$.

In the situation of theorem \ref{Tm:Coprod} it makes sense to ask
if there exists an antipode
for arbitrary choice of the coefficients $q_{ij}$. An explicit
description of $S$ follows by a general argument. A Hopf
\textsf{Ass}-algebra $A$ is called \note{connected} if it is $\ZZ$-graded,
concentrated in non-negative degree, and satisfies $A_0=k\cdot 1$. The
\note{augmentation ideal} of a connected Hopf algebra $A$ is the ideal
$\bigoplus_{n\geq 1} A_n$.

\begin{Lm}[Milnor and Moore \cite{MilnorMoore:Hoa}]
Let $A$ be a connected Hopf \textsf{Ass}-algebra. Then there exists an
antipode on $A$ which for $x$ in the augmentation ideal is given by
\[
S(x)= \sum_{k=1}^\infty (-1)^{k+1}\mu^{(k)}\circ\bar{\Delta}^{(k)}(x),
\]
where $\bar{\Delta}= \Delta - (\id\otimes 1 + 1 \otimes \id)$, and
$\mu^{(0)}=\id = \bar{\Delta}^{(0)}$, and $\mu^{(k)}:A^{\otimes k+1}\rightarrow A$ and
$\bar{\Delta}^{(k)}:A\rightarrow A^{\otimes k+1}$ are defined using
(co)associativity for $k>0$ and $\mu^{(0)}=\id = \bar{\Delta}^{(0)}$.
\end{Lm}


\begin{Cr}
Under the assumptions of convention \ref{Cv:Deform}, there exists an antipode
on the Hopf algebra of theorem \ref{Tm:Coprod}, given by the formula 
\[
S(t) = \sum _{k=1}^{|t|}\sum_{\cup_i s_i = t} (-1)^k s_1\clcdot
s_k 
\prod_{1\leq j <k}q(s_j,s_{j+1}\cup\ldots\cup s_k,s_j\cup\ldots\cup s_k),
\]
where we only sum over partitions $t = s_1\cup\ldots \cup s_k$ of the
forest $t$ with all forests $s_i$ non-empty.
\end{Cr}
\begin{Pf}
The grading with respect to $|t|$ gives a grading that makes the Lie algebra 
connected, since $\textsf{Com}(0)=k$.
\end{Pf}
\begin{Ex}
In the case $n=1$ with $q_1=1$ and $q_2 = 0$ we get the formula of
Connes and Kreimer \cite{ConKr:Hoa} for the antipode on the Hopf
algebra $H_R$ of section \ref{Sec:Preliminaries}. Note that according
to example \ref{Ex:Twistfunc} the Connes-Kreimer Hopf algebra is the
initial Hopf algebra for the Hopf operad $\textsf{Com}[\lambda_1]$.
\end{Ex}

We now state a corollary to the formula for the coproduct in theorem
\ref{Tm:Coprod}, that gives an explicit formula for the Lie bracket on
the primitive elements of the cocommutative Hopf algebra $C_n^*$. We
know by the Milnor-Moore theorem \cite {MilnorMoore:Hoa} that $C_n^*$ is
the universal enveloping algebra of the Lie algebra of its primitive
elements.

\begin{Cr}\label{HnLiefomrula}
Let $\textsf{P} = \textsf{Com}$, and $t\in C_n$, and
$\Delta = \Delta_{q_{1\,1}\ldots q_{2\,n}}$; then the graded dual $C_n^*$
is the universal enveloping algebra of the Lie algebra which as a
vector space is spanned by elements $D_t$, where $t$ is a rooted tree in $C_n$.
The bracket is given by $[D_s,D_t] = D_s*D_t-D_t*D_s$, where
\[
D_s*D_t = \sum_{w=s\cup t} q(s,w)\, D_w,
\]
the sum ranges over all rooted trees in $C_n$ that have $s$ 
and $t$ as complementary colored subforests.
\end{Cr}

\begin{Pf}
For any cocommutative Hopf algebra we can define an operation $*$ on
the primitive 
elements, such that its commutator is the Lie bracket on primitive
elements: Simply 
define $*$ as the truncation of the product at degree $>1$, with
respect to the primitive 
filtration $F$. In this case, $F_mC^*_n$ is spanned by the elements
$D_u$ dual to forests $u$ consisting of at most $m$ trees.
The product in $C_n^*$ is determined by the coproduct in
$C_n$. Explicitly, for forest $u$, we can write the multiplication in
$C^*_n$ as
\[
\begin{split}
D_s D_t (u) &= (D_s \otimes D_t) \Delta(u) \\
&= \sum_{u=w\cup w^c} q(w,u) D_s(w)D_t(w^c).
\end{split}
\]
When we then restrict to the primitive part, we conclude that
$D_s*D_t$ is given by the desired formula.
\end{Pf}

\begin{Rm}\label{Rm:Assdeform}
With some minor changes, there is an analogue of the above for the case 
\textsf{P = Ass}. The only change is that we have to remember the ordering 
of up-going edges at each vertex. To be very explicit, for the coproduct on 
$A_n = \textsf{Ass}[\lambda_n](0)$ we have the formula
\[
\Delta(t) = \sum_{s\subset t} 
\prod_j q_{1\,j}^{\sum_{v\in s}p_j(v,s,t)}\cdot
\prod_jq_{2\,j}^{\sum_{v\in t}p_j(v,s^c,t)}
s\otimes s^c,
\]
where the product of trees is associative, non-commutative. The order of
multiplication is given by the linear order on the roots of the trees 
defined by the linear on the incoming edges at each vertex and the partial 
order on vertices. Dually, the primitive elements of the dual are again 
spanned by elements dual to single rooted trees (with linear ordering
on incoming edges in each vertex). The truncated product 
$*$ is given by
\[
D_s*D_t = \sum_{w=s\cup t} 
\prod_j q_{1\,j}^{\sum_{v\in s}p_j(v,s,w)}\cdot
\prod_jq_{2\,j}^{\sum_{v\in t}p_j(v,t,w)}
D_w,
\]
where $w$, $s$ and $t$ are trees with a linear ordering on the incoming 
edges of the same color at each vertex and the inclusions of $s$ and $t$ in 
$w$ have to respect these orderings.
\end{Rm}

\begin{Ex}
Independently, Foissy \cite{Fois:Plantree} has found this formula 
for the Lie bracket in  the case where  $n=1,\ q_1=1$ and $q_2=0$ (and 
$\textsf{Ass}$ is the underlying operad), and established the universal 
property of section \ref{Sec:Cohomology}, which in this case also follows 
from Moerdijk \cite{Moer:hoa}. He uses this formula to give an explicit 
isomorphism between the Hopf algebras $A_1$ and $A_1^*$ in this case. 
Moreover, he does the same for planar trees with labeled vertices. The
labeled case 
corresponds to iterated addition of unitary operations, using the same 
$q_1$ and $q_2$ at each stage.
\end{Ex}

\begin{Rm}\label{Rm:kdeform}
For certain values of the parameters $q_{ij}$, the formula of theorem \ref{Tm:Coprod} defines
a coproduct on $\textsf{k}[\lambda_n](0)$.
If for fixed $i=1,2$ there is at most one $j$ such that $q_{ij}\neq 0$, then 
the subspace spanned by single rooted trees which do not involve 
multiplication is preserved under $\Delta$. Since this subspace is 
isomorphic to $\textsf{k}[\lambda_n](0)$, this implies that in this case the 
formulas above define a coproduct on $\textsf{k}[\lambda_n](0)$. 
Consequently, this defines an operation $*$ on the dual, since each element 
in the dual is primitive. A formula for this is given by the formula for 
\textsf{Com}, modulo trees involving multiplication (their span is denoted $\mu$):
the quotient $\textsf{Com}[\lambda_n](0)^*/(\mu)$ is canonically isomorphic 
to $\textsf{k}[\lambda_n](0)^*$.
\end{Rm}
%
%

\section{Pre-Lie algebras}
\label{Sec:PreLie}

A \note{(right) pre-Lie algebra} (cf. Gerstenhaber \cite{Gerst:Ring}, and Chapoton and
Livernet \cite{ChapLiv:Prelie}) 
is a $k$-vector space $L$ together with a bilinear map 
$-*-: L \times L \rarr L$ which satisfies
\[
(x*y)*z - x*(y*z) = (x*z)*y - x*(z*y),
\]
for each $x,\ y$ and $z$ in $L$. Note that the bracket
$[x,y]=x*y-y*x$ defines a Lie algebra structure on $L$. This Lie
algebra is the \note{associated Lie Algebra} of the pre-Lie
algebra $(L,*)$.

Remember from \cite{ChapLiv:Prelie} that the free pre-Lie algebra
$L_p$ on $p$ generators is given by the vector space spanned by rooted trees 
with vertices labeled with elements of $\mathbf{p}=\{1,2,\ldots,p\}$. 
The pre-Lie algebra product is given by grafting trees. That is,
\[
t*s = \sum_{v\in t} s\circ_v t
\]
for trees $s$ and $t$, where the color of the vertices is
preserved. For any tree $t\in L_p$ define $r_i(t)$ to be the tree
$t$ with the color of the root changed to $i$. Thus $r_i$ defines
a linear endomorphism of the vector space $L_p$.

\begin{Cv}
Let \textsf{P = Com}, choose $\Delta$ on $C_n = \textsf{Com}[\lambda_n](0)$ 
defined by $\sigma_i$ as in the deformation used in convention \ref{Cv:Deform} 
and use the notation $\chi_S$ to denote the characteristic function of a
subset $S\subset X$ which has value $1$ on $S$ and value $0$ on
$X-S$. Denote the primitive elements of a coalgebra $C$ by $P(C)$.
\end{Cv}
Let $t \in C_n$ be a tree. Then we call the $\mathbf{p}$-connected component 
of $t$ the 
subtree consisting of vertices connected to the root by a path with edges 
colored exclusively by elements of $\mathbf{p}$.

\begin{Pp}\label{Tm:preLie}
Let $\mathbf{p}\subset\{1,\ldots, n\}$ and define
$q_{1j}=\chi_{\mathbf{p}}(j)$ and $q_{2j}=0$. Then
\begin{enumerate}
\item The product $*$ of corollary \ref{HnLiefomrula} defines a pre-Lie
algebra structure on the primitive elements $P(C^*_n)$ of $C_n^*$.
\item If $\mathbf{p}=\{1,\ldots,n\}$, then there is a natural inclusion
of this pre-Lie algebra into the free pre-Lie algebra on $n$
generators. The image in the free pre-Lie algebra is spanned by
all sums $\sum_{i\in\mathbf{p}}t_i$, of trees that only differ in
that the color of the root of $t_i$ is $i$.
\end{enumerate}
\end{Pp}

\begin{Pf}
First note that with this choice of $\sigma_i$, the formula for
the coproduct simplifies to
\[
\Delta(t) = \sum_{c\in C(t)} P^c(t)\otimes R^c(t),
\]
where we sum over admissible cuts as described in section
\ref{Sec:Preliminaries},
provided that the cuts are at edges in the $\mathbf{p}$-connected component of $t$.
To see this, consider the formula of corollary
\ref{HnLiefomrula} and note that the coefficient of $s\otimes s^c$ is 0 iff
$s^c$ has a connected component with respect to edges of colors in 
$\mathbf{p}$ which does not contain the root. 

The operation $*$ of theorem
\ref{HnLiefomrula} is then given by
\[
D_s*D_t = \sum_{u} \sum_{p\in\mathbf{p}} n_p(s,t,u) D_{u},
\]
where  $u$ ranges over trees, and $n_p(s,t,u)$ 
denotes the number of admissible cuts $c$ of $u$ that only cut  $p$-colored 
edges in the $\mathbf{p}$-connected component of $u$ and satisfy 
$P^c(u) = t$ and $R^c(u)=s$.

To check the
pre-Lie algebra identity, note that $D_u*(D_t*D_s) - (D_u*D_t)*D_s$ 
is exactly  $\sum_w n(u,t,s,w) D_w$, where $n(u,t,s,w)$ is the number of 
admissible cuts $c$ of $w$, such that $P^c(w) = u \cdot t$ and $R^c(w) = s$.
Since $n(t,u,s,w) = n(u,t,s,w)$ part \textit{(i)} follows.

For a tree $t$, we denote the cardinality of the automorphism 
group of $t$ by $a_t$. If $D_t \in P(C_n^*)$ and $v$ a vertex of $t$, 
denote by $a_t^v$ the cardinality of the  group of automorphisms of $t$ 
that have $v$ as a fixed point. Define $m(v,t)$ the cardinality of the 
orbit of $v$ under the action by the automorphisms of $t$. Let $s, t\in C_n$ 
be trees and define $u = s \circ_v^p t$. We now see that
\[
a_u = \sum_p a_s \cdot a_t^v \cdot n_p(s,t,u).
\]
Let $\mathbf{p}=\{1,\ldots,n)$ and consider the pre-Lie algebra of
\textit{(i)}. We claim that there is an isomorphism of pre-Lie algebras from 
this pre-Lie algebra to the pre Lie algebra on the same vector space $C_n$ 
with multiplication given by
\[
D_s*'D_t = \sum_{v\in t}\sum_p D_{s\circ_v^p t}.
\] 
Define a linear automorphism of $P(C_n^*)$ by $\phi(D_t) = a_tD_t$.
We have to show that $\phi(D_s*D_t) = \phi(D_s)*'\phi(D_t)$. Write
\[
\begin{split}
D_s*'D_t &= \sum_{v\in t}\sum_p D_{s\circ_v^p t}\\ 
                &= \sum_p\sum_{u=s\circ_v^p t}D_u m(v,t)
\end{split}
\]  
where the sum in the second line is over all $u$ that can be written
as $s\circ_v^p t$ for some $v$ and where for each such $u$ we 
fix a decomposition of each such $u$ as $u = s \circ_v^p t$. To show that 
$\phi$ is a pre-Lie algebra morphism, we have to show that for any such $u$,
\[
\sum_p a_s \cdot a_t \cdot n_p(s,t,u) = a_u m(v,t).
\] 
This follows from the above since $a_t = a_t^v \cdot m(v,t)$.

 A map from the pre-Lie algebra $(P(C_n^*),*')$ into the free pre-Lie algebra
on $n$ generators (trees with $n$-colored vertices) is given by
assigning the color of an edge to the vertex on top and summing
over all colors in $\mathbf{p}$ for the root. That this is an inclusion 
of pre-Lie algebras and that the image can be described as in \textit{(ii)} is a 
direct consequence of the definitions.
\end{Pf}

\begin{Ex}
In the case $n=1$, $\mathbf{p}=\{1\}$, this yields the identification
of $L_1$ with the primitive elements of $H_R^*$ with
the pre-Lie product $*$ as in Chapoton and Livernet \cite{ChapLiv:Prelie}. 
\end{Ex}


%
%
\section{The Pruning coproduct}
\label{Sec:Pruning}
Brouder and Frabetti \cite{BrouFrab:Trees} use a pruning operator on
planar binary trees to define a 
coproduct on the $k$-vector space spanned by planar binary trees. We
will see that this coproduct is 
related to the coproducts described in section \ref{Sec:Deformation}. We
first give a family of coproducts on
planar $n$-ary trees. In this context the pruning Hopf algebra of
Brouder and Frabetti is a direct corollary.

\label{Txt:bintrees}
A \note{planar binary tree} is a non-empty oriented planar graph in
which each vertex which is not a leaf 
or the root has exactly two direct predecessors and one direct
successor. Let $T$ be a planar binary 
tree. We call each vertex which is not a leaf a \note{internal vertex}
of $T$. If $T$ has $n+1$ leaves, it has $n$ internal
vertices. Following Loday and Ronco, we denote the $k$-linear span of
the set of planar binary trees $k[Y_\infty]$. On planar binary trees we have
the binary operation $\lambda$ adding a new root (smallest element) to
the union of planar binary trees $T_1$ and $T_2$. 

Likewise, we define a \note{ planar $n$-ary tree} as a non-empty
oriented graph in which each vertex which is not a leaf has exactly
$n$ direct predecessors and one direct successor. The $n$-ary
operation $\lambda$ adds a new root to an $n$-tuple
of planar $n$-ary trees. Denote the vector space spanned by the set of
$n$-ary planar trees by $k[Y^{(n)}_\infty]$. The basepoint in
$k[Y^{(n)}_\infty]$ is the tree with one leaf.

\begin{Pp}\label{Pp:freeass}
Let $n\geq 0$. 
\begin{enumerate}
\item 
The pointed vector space $k[Y^{(n)}_\infty]$ is isomorphic 
to $\textsf{k}[\lambda_n](0)$.
\item
The unitary free associative algebra on the pointed vector space spanned by
planar $n$-ary trees has a family of natural Hopf algebra structures
parameterized by $i\in\{1,\ldots,n\}$.
\end{enumerate}
\end{Pp}
\begin{Pf}
Define a bijection $\xi: k[Y^{(n)}_\infty]\rightarrow
\textsf{k}[\lambda_n](0)$ by mapping the planar $n$-ary tree with one leaf 
to the empty tree and 
\[
\xi(\lambda(T_1,\ldots,T_n)) = \lambda_n(\xi(T_1),\ldots,\xi(T_n)),
\]
for planar $n$-ary trees $T_i$. Each vertex in a tree in
$\textsf{k}[\lambda_n](0)$ corresponds to an 
internal vertex in the corresponding $n$-ary tree.
It is easy to check that this is a basepoint preserving bijection of bases.
This proves \textit{(i)}.

Let $i\in\{1,\ldots,n\}$. The coproducts induces by $q_{1j}=\chi_{\{i\}}(j)$ 
and $q_{2j}=0$ are well defined on $\textsf{k}[\lambda_n](0)$ by
remark \ref{Rm:kdeform}.
This shows that we have a family of coproducts on $\textsf{k}[\lambda_n](0)$.
The result \textit{(ii)} is now a direct consequence of the
deformation in remark \ref{Rm:Assdeform} and functorial properties.
The functor 
$\bar{i}_!:\textsf{k}[\lambda_n]\text{-HopfAlg}
\rightarrow\textsf{Ass}[\lambda_n]\text{-HopfAlg}$ 
induced by the inclusion of Hopf operads 
$i:\textsf{k}\rightarrow \textsf{Ass}$
is the unitary free algebra functor and the result follows by 
functoriality.
\end{Pf}

\begin{Cr}[Brouder and Frabetti \cite{BrouFrab:Trees}]\label{Pp:BrouFrab}
Let the pruning operator $P:k[Y_\infty]\rarr k[Y_\infty]\otimes k[Y_\infty]$ be defined
inductively as
\[
\begin{split}
P(\lambda(T,1))&=0 \\
P(\lambda(T,S) &= \sum_i\lambda(T,S_i')\otimes S_i'' + \lambda(T,1)\otimes S,
\end{split}
\]
if $P(S) = \sum_i S_i'\otimes S_i''$. Then
the formulas
\[
\begin{split}
\Delta^P(1) &= 1\otimes 1 \quad\text{and}\\
\Delta^P(T) &= 1\otimes t + P(T) + T\otimes 1 ,
\end{split}
\]
define a coassociative coproduct on $k[Y_\infty]$. Thus the
unitary free associative algebra $T(k[Y_{\infty}])$ is a Hopf algebra
with the 
coproduct described on generators by $\Delta^P$.
\end{Cr}
\begin{Pf}
Let $n=2$, and $i=2$. Apply proposition \ref{Pp:freeass}.
\end{Pf}

\section{Associative and dendriform algebras}

\label{Sec:Associative}
A \note{dendriform algebra} (cf. Loday \cite{Loday:bialgebras}, and
Loday and Ronco
\cite{LodRon:Trees}) is a vector space together with two bilinear
(non-associative) products $\prec$ and $\succ$, satisfying the
identities
\[
\begin{split}
(a\prec b) \prec c &= a\prec (b\prec c) + a \prec (b \succ c) \\
a \succ (b \prec c) &= (a \succ b) \prec c \\
a\succ (b \succ c) &= (a\prec b) \succ c) + (a \succ b) \succ c.
\end{split}
\]
A dendriform algebra $D$ defines an associative algebra on the
same vector space with associative product $*$ defined by
\[
a * b = a \prec b + a \succ b.
\]
A \note{dendriform Hopf algebra} (Ronco \cite{Ron:MiMo}) is a dendriform
algebra together with a coproduct that satisfies 
\[
\begin{split}
\Delta(x\prec y) &=\sum x' * y' \otimes x''\prec y'' + x'*y\otimes x''
+ y'\otimes x\prec y'' + y\otimes x \\
\Delta(x\succ y) &=\sum x'*y'\otimes x''\succ y'' + x*y'\otimes y'' +
y'\otimes x\succ y'' + x\otimes y. 
\end{split}
\]
\begin{Rm}
One motivation for studying dendriform algebras is the
classification of Hopf algebras analogous to the Milnor-Moore
theorem (cf. Milnor and Moore \cite{MilnorMoore:Hoa}). The
Milnor-Moore theorem gives an equivalence of categories between
the category of cocommutative Hopf algebras and the category of
Lie algebras (in characteristic 0). A Milnor-Moore theorem for
dendriform Hopf algebras can be found in Ronco \cite{Ron:MiMo}. This
theorem states the equivalence of dendriform Hopf algebras and brace
algebras.
\end{Rm}

\begin{Cv}\label{Cv:LodRon}
Let $k,l\in\NN$. Use \textsf{P = k}, and take
$q_{1j}=\delta_{kj}$ and $q_{2j}=\delta_{lj}$ to define $\sigma_1$
and $\sigma_2$ in the deformation of \ref{Tm:Coprod}. 
According to remark \ref{Rm:kdeform} this choice of parameters indeed
defines a coproduct on $\textsf{k}[\lambda_p](0)$. 
\end{Cv}
We use the augmentation of $u^*:\textsf{k}[\lambda](0)^*\rightarrow k$.
\begin{Pp}\label{Pp:Starass}
Consider the graded vector space $(\textsf{k}[\lambda_p](0))^*$, graded by the 
number of applications of $\lambda_p$ in trees. Then 
the augmentation ideal of $(\textsf{k}[\lambda_p](0))^*$ is a graded dendriform
algebra with associated multiplication $*$.
\end{Pp}

\begin{Pf}
The most convenient way to check the statement is to give an
inductive characterization of the 
operations $\prec$ and $\succ$. To avoid awesome
notation, we identify trees $t$ with their dual $D_t$.

The product is characterized by its unit $1=\emptyset\
(=D_{\emptyset})$, and the identity
\[
\begin{split}
s*t &= \lambda(s_1,\ldots,s_k*t,\ldots,s_n) +
\lambda(t_1,\ldots,s*t_l,\ldots,t_n),
\end{split}
\]
as is readily checked. Now define operations
\[
\begin{split}
s \prec t &= \lambda(s_1,\ldots,s_k*t,\ldots,s_n) \quad \text{and}\\
s \succ t &= \lambda(t_1,\ldots,s*t_l,\ldots,t_n).
\end{split}
\]
Let $t=\lambda(t_1,\ldots,t_n)$, $u=\lambda(u_1,\ldots,u_n)$ and
$s=\lambda(s_1,\ldots,s_n)$ and suppose that $\prec$ and $\succ$
satisfy the dendriform identities on all trees with in total $<|t|
+ |s| +|u|$ vertices. Then we have
\[
\begin{split}
(s\prec t) \prec u &= \lambda(s_1,\ldots,(s_k*t)*u,\ldots,s_n)\\
&=\lambda(s_1,\ldots,s_k*(t*u),\ldots,s_n)\\
&=s\prec(t*u) \\
(s\succ t) \prec u &= \lambda(t_1,\ldots,s*t_l,\ldots,
t_k*u,\ldots,t_n)\\
&= s\succ(t\prec u).
\end{split}
\]
The first identity follows by the definition of $*$ and
associativity on smaller trees. The second is clear, but uses
associativity on smaller trees in case $k=l$. The proof of the third 
dendriform identity can be copied from the first almost verbatim. 
It is now left to show that the dendriform operations preserve the grading
in the sense that for $s$ of degree $m$ and $t$ of degree $n$, the
products $s\prec t$ and $s\succ t$ are of degree $m+n$. If the
dendriform operations are graded in total degree smaller then
$m+n$, assume $s$ and $t$ of degree $m$ and $n$ respectively. We
now know that $*$ is graded in total degree smaller then $m+n$.
Since the definition of $\prec$ and $\succ$ only uses
multiplication of smaller trees, one checks that the grading is preserved.
\end{Pf}

\section{The Loday-Ronco construction}

\label{Sec:LodayRonco}

We still assume convention \ref{Cv:LodRon}.
There is a isomorphism of graded vector spaces between
$\textsf{k}[\lambda_n](0)$ and $(\textsf{k}[\lambda_n](0))^*$ that
maps $t$ to $D_t$. Proposition \ref{Pp:Starass} shows that we have an
associative multiplication $*$ on
$(\textsf{k}[\lambda_n](0))^*$. The isomorphism give an associative
multiplication on $\textsf{k}[\lambda_n](0)$, which we will also
denote by $*$.

Denote by $*_{(n)}:\textsf{k}[\lambda_n](0)^{\otimes n} \rightarrow
\textsf{k}[\lambda_n](0)$ the map induced by the associative
multiplication, and by $\eps_{(n)}:\textsf{k}[\lambda_n](0)^{\otimes
n} \rightarrow \textsf{k}[\lambda_n](0)$ the map $u\circ\eps^{\otimes n}$.

\begin{Pp}\label{Pp:LRgen}
Use convention \ref{Cv:LodRon} and let the product $*$ of proposition 
\ref{Pp:Starass}. Then
$\sigma_1=*_{(n)}$ and $\sigma_2 = \eps_{(n)}$ define a Hopf algebra
structure on $\textsf{k}[\lambda_n](0)$ iff $k=n$ and $l=1$.

Moreover if $k=n$ and $l=1$, the augmentation ideal of
$\textsf{k}[\lambda_n](0)$ is a dendriform Hopf algebra.
\end{Pp}
\begin{Pf}
We check the requirements of proposition \ref{Tm:genmoerdijk}. The
requirement for the counit is fulfilled for any element of
$(\textsf{k}[\lambda_p](0))^*$. We prove the claim for coassociativity
inductively. Let $\sigma_i$ as above and assume that the requirement
for $\Delta$ holds up to degree $n$. Then
coassociativity of $\Delta$ holds up to degree $n$. Compute for $k=n$
and $l=1$ (using
Sweedler notation)
\[
\begin{split}
\Delta(a*b) &= 
\Delta(\lambda(a_1,\ldots,a_n)*\lambda(b_1,\ldots,b_n)) \\
&= \Delta(\lambda(a_1,\ldots,a_n*b) 
+ \Delta(\lambda(a*b_1,\ldots,b_n))\\
&= \sum
a'*\ldots *(a_n*b)' \otimes\lambda(a_1'',\ldots,(a_n*b)'')\\
& + \sum
(a*b_1)'*\ldots *b_n'\otimes\lambda((a*b_1)'',\ldots,b_n'')\\
&+ (a*b)\otimes 1
\end{split}
\]
for $a=\lambda(a_1,\ldots,a_n)$ and $b=\lambda(b_1,\ldots,b_n)$ such that 
the sum of their degrees is $m$. Then use 
the induction hypothesis for degrees smaller then $m$
(in particular $\Delta(a_n*b) = \Delta(a_n)*\Delta(b)$
and $\Delta(a*b_1) = \Delta(a)*\Delta(b_1)$). Deduce
\[
\begin{split}
\Delta(a\prec b) &= \sum a'*b'\otimes a''\prec b''\qquad \text{ and }\\ 
\Delta(a\succ b) &= a'*b'\otimes a''\succ b'',
\end{split}
\]
and conclude that
$(\textsf{k}[\lambda_p](0))^*$ is a Hopf algebra if $k=n$ and $l=1$.
The formulas also show that for other choices of $k$
and $l$ the induction fails and coassociativity does not hold since
$*$ is not commutative. 

The dendriform Hopf identities for
$\bar{\Delta} = \Delta - \id\otimes 1 - 1 \otimes \id$ follow from the 
above.
\end{Pf}

Loday and Ronco \cite{LodRon:Trees} construct a
Hopf algebra of planar binary trees. Remember the definition of
the vector space $k[Y_{\infty}]$ as the vector space spanned 
by binary trees.

\begin{Cr}[Loday and Ronco \cite{LodRon:Trees}]
Let $T , T' \in k[Y_\infty ]$. If $T=\lambda(T_1,T_2)$ and
$S=\lambda(S_1,S_2)$, then
\[
T*S = \lambda(T_1,T_2*S) + \lambda(T*S_1S_2)
\]
defines an associative product on $k[Y_{\infty}]$ with as unit the
tree with one leaf. 

Moreover there is a coproduct making
$k[Y_{\infty}]$ a Hopf algebra. This coproduct is defined
inductively by
\[
\Delta(T) = \sum (T_{1}'*T_{2}')\otimes
\lambda(T_{1}'',T_{2}'') + T \otimes 1,
\]
where $\Delta(T) = \sum_{(T)} T'\otimes T''$ is the
Sweedler notation of the coproduct. The counit $\eps$ is the
obvious augmentation of the product by projection on the span of
the unit.

Finally, the augmentation ideal is a dendriform 
Hopf algebra. 
\end{Cr}
\begin{Pf}
Recognize the star product on
$\textsf{k}[\lambda_2](0)$ given by $q_{12}=1=q_{21}$ and
$q_{11}=0=q_{22}$, and the coproduct associated to $*$. Apply
proposition \ref{Pp:LRgen}.
\end{Pf}

\bibliographystyle{plain}
\bibliography{hopf}

\begin{thebibliography}{10}

\bibitem{BrouFrab:Trees}
C.~Brouder and A.~Frabetti.
\newblock Renormalization of {QED} with {P}lanar {B}inary {T}rees.
\newblock {\em Euro. Phys. J. C}, 19(4):715--741, 2001.

\bibitem{ChapLiv:Prelie}
F.~Chapoton and M.~Livernet.
\newblock Lie {A}lgebras and the {R}ooted {T}rees {O}perad.
\newblock {\texttt{arXiv:math.QA/0002069}}, 2000.

\bibitem{ConKr:Hoa}
A.~Connes and D.~Kreimer.
\newblock Hopf {A}lgebras, {R}enormalisation and {N}oncommutative {G}eometry.
\newblock {\em Comm. Math. Phys.}, 199, 1998.

\bibitem{Fois:Plantree}
L.~Foissy.
\newblock Les {A}lg\`ebres de {H}opf des {A}rbres {E}nracin\'es {D}ecor\'es.
\newblock \texttt{arXiv:math.QA/0105212}, 2001.

\bibitem{Gerst:Ring}
M.~Gerstenhaber.
\newblock The {C}ohomology {S}tructure of an {A}ssociative {R}ing.
\newblock {\em Ann. Math.}, 78(2):59--103, 1963.

\bibitem{GetzJon:Opd}
E.~Getzler and J.D.S. Jones.
\newblock Operads, {H}omotopy {A}lgebra, and {I}terated {I}ntegrals for
  {D}ouble {L}oop {S}paces.
\newblock {\texttt{arXiv:hep-th/9403055}}, 1994.

\bibitem{GinKap:Koszul}
V.A. Ginzburg and M.M. Kapranov.
\newblock Koszul {D}uality for {O}perads.
\newblock {\em Duke Math. J.}, 76:203--272, 1994.

\bibitem{KrizMay:Opd}
I.~Kriz and J.~May.
\newblock Operads, {A}lgebras, {M}odules and {M}otives.
\newblock {\em Ast\'erisque}, 233, 1995.

\bibitem{LazMov:hoa}
A.Yu. Lazarev and M.V. Movshev.
\newblock Deformations of {H}opf {A}lgebras.
\newblock {\em Russian Math. Surveys}, 46(1):253--254, 1991.

\bibitem{Loday:bialgebras}
J.L. Loday.
\newblock Dialgebras.
\newblock {\texttt{arXiv:math.QA/0102053}}, 2001.

\bibitem{LodRon:Trees}
J.L. Loday and Maria~O. Ronco.
\newblock Hopf {A}lgebra of the {P}lanar {B}inary {T}rees.
\newblock {\em Adv. Math.}, 139(2):293--309, 1998.

\bibitem{MilnorMoore:Hoa}
J.~Milnor and J.~Moore.
\newblock On the {S}tructure of {H}opf {A}lgebras.
\newblock {\em Ann. Math}, 81(2):211--264, 1965.

\bibitem{Moer:hoa}
I.~Moerdijk.
\newblock On the {C}onnes-{K}reimer {C}onstruction of {H}opf {A}lgebras.
\newblock {\em Contemp. Math.}, 271:311--321, 2001.
\newblock Also at {\texttt{arXiv:math-ph/9907010}}.

\bibitem{Ron:MiMo}
M.~O. Ronco.
\newblock A {M}ilnor-{M}oore {T}heorem for {D}endriform {H}opf {A}lgebras.
\newblock {\em C.R. Acad. Sci. Paris S\'er.I Math.}, 332(2):109--114, 2001.

\bibitem{Sweedler:Hoa}
M.E. Sweedler.
\newblock {\em Hopf {A}lgebras}.
\newblock Oxford University Press, 1969.

\end{thebibliography}

\par
Pepijn van der Laan \hfill P.O.Box 80.010\\
Mathematisch Instituut \hfill 3508 TA Utrecht\\
Universiteit Utrecht \hfill The Netherlands\\
\end{document}